\documentclass[12pt, a4paper]{article}
\usepackage[T1]{fontenc}
\usepackage{setspace}
\doublespace
\begin{document}
\renewcommand{\baselinestretch}{1.5}
\newtheorem{Definition}{Definition}[section]
\newtheorem{Theorem}{Theorem}[section]
\newtheorem{Lemma}[Theorem]{Lemma}
\newtheorem{Corollary}[Theorem]{Corollary}
\title{Random Interval Graphs}
\author{Iliopoulos Vasileios}
\date{}
\maketitle
\tableofcontents
\pagebreak

\noindent
\section{Abstract}
In this thesis, which is supervised by Dr. David Penman, we examine
random interval graphs. Recall that such a graph is defined by
letting $X_{1},\ldots X_{n},Y_{1},\ldots Y_{n}$ be $2n$ independent
random variables, with uniform distribution on $[0,1]$. We then say
that the $i$th of the $n$ vertices is the interval $[X_{i},Y_{i}]$
if $X_{i}<Y_{i}$ and the interval $[Y_{i},X_{i}]$ if $Y_{i}<X_{i}$.
We then say that two vertices are adjacent if and only if the
corresponding intervals intersect.

We recall from our MA902 essay that fact that in such a graph, each
edge arises with probability $2/3$, and use this fact to obtain
estimates of the number of edges. Next, we turn to how these edges
are spread out, seeing that (for example) the range of degrees for the
vertices is much larger than classically, by use of an interesting
geometrical lemma. We further investigate the maximum degree, showing
it is always very close to the maximum possible value $(n-1)$, and the
striking result that it is equal to $(n-1)$ with probability exactly $2/3$.
We also recall a result on the minimum degree, and contrast all these
results with the much narrower range of values obtained in the alternative
\lq comparable\rq\, model $G(n,2/3)$ (defined later).

We then study
clique numbers, chromatic numbers and independence numbers in the
Random Interval Graphs, presenting (for example) a result on independence
numbers which is proved by considering the largest chain in the associated
interval order.

Last, we make some brief remarks about other ways to define random
interval graphs, and extensions of random interval graphs, including
random dot product graphs and other ways to define random interval
graphs. We also discuss some areas these ideas should be usable in. We close with a summary and
some comments.
\pagebreak

\noindent
\section{Acknowledgements}
I would like to thank my supervisor, Dr. David Penman for his
motivation and support. Also, I would like to express my sincere
thanks to all my teachers in this academic year. All mistakes and
errors in this work are my own.
\pagebreak

\noindent
\section{Basic results}
\subsection{Definitions}
We recall here the definitions of interval graph, and random interval
graph, from our essay \cite{iliop}.
\begin{Definition}
Let a graph $G$ have $n$ vertices $\{1,2,\ldots n\}$. To create an Interval
Graph, to each vertex $i$ we assign a finite interval $I_{i}$of the real line.
This transformation yields $n$ intervals of the real line. We then say that
two different vertices of $G$ are adjacent, in the initial graph, if the
corresponding intervals have non-empty intersection. That is,
$i\sim j\leftrightarrow I_{i}\cap I_{j}\neq \emptyset$.
\end{Definition}
\begin{Definition}
A Random Interval Graph is formed as follows. Suppose the vertex set is
$\{1,2\ldots n\}$. Let $2n$ independent and identically distributed
continuous random variables, $X_{1}, X_{2}\ldots X_{n}$ and $Y_{1},
Y_{2}\ldots Y_{n}$, from the uniform distribution in $[0, 1]$, be given. The
interval $I_{i}$ will be $[X_{i}, Y_{i}]$, if $X_{i}<Y_{i}$ or
$[Y_{i}, X_{i}]$, when $X_{i}>Y_{i}$. (The case where any two random
variables are equal has probability 0, by properties of continuous
distributions). We then say that the random interval graph is the interval
graph formed for these vertices from the intervals $I_{i}$.

We use $\Delta_{n}$ to denote the set of all possible random interval
graphs on $n$ vertices.
\end{Definition}
The main aim of this essay will be to prove various basic properties
of these graphs, based on the two papers \cite{Schein} by E. R.
Scheinerman and \cite{JSTOR}. We will first show some equivalent
formulations of the model, which will be useful in proving theorems,
and will then prove some results about various aspects of the
graphs.

\subsection{An Equivalent model}
In \cite{Schein} it is observed that we do not have to use the
uniform random variables $X_{i}$ and $Y_{i}$ above. He observed that
it is enough to suppose that the intervals have as their endpoints
the numbers 1,2\ldots 2n in some random order, with all the $(2n)!$
possible orderings equally likely. The reason why this is equivalent
is that, as we observed when defining random interval graphs, the
probability that two of the $X_{i}$ and $Y_{j}$ are equal is zero,
so they can be any set of unequal numbers. Also, because the $X_{i}$
and $Y_{i}$ are all independent, all $(2n)!$ possible orderings of
them have the same probability. (Independence implies
exchangeability - that is, the property that the probability that
the $X_{i}$s and $Y_{i}$s take certain values is the same as the
probability that the image of them all under some permutation in
$S_{2n}$, the symmetric group on $2n$ letters, take these values).
This concludes a proof of the following result, which will be used
later in the proof that there is vertex of degree $(n-1)$ with
probability $2/3$.
\begin{Lemma}
An equivalent definition of random interval graphs is to say they
have vertex set $\{1,2\ldots n\}$ and that the intervals attached to
the vertices have as their $2n$ endpoints the numbers $1,2\ldots 2n$
in some order, the order being chosen uniformly at random so that
all $(2n)!$ possible orders are equally likely.
\end{Lemma}
{\bf Proof.} See above.

\subsection{Edges in interval graphs}
The simplest question about a graph is how many edges it has.
By Theorem $10$ in \cite{iliop}, we know that the probability of any
particular edge arising in a random interval graph is $2/3$. We first
extend this result to obtain an estimate of the total number of edges
in a random interval graph, namely that it is near to $n^{2}(1+o(1))/3$.
This is the first main result in E. R. Scheinerman's paper \cite{Schein}.

It will be helpful to give an overview of the proof. It is easy to show that
the expected number of edges is exactly $n(n-1)/3$ using the fact about the
probability that an individual edge arises is $2/3$. What we need to do is
to show that it is very likely to be very close to this value - the precise
sense of this will be made clear in the statement of the result. (The result
is a limit result: many interesting results in probability, e.g. the law of
large numbers and the Central Limit Theorem, are of this form). In order
to show that it is very likely to be very close to $n(n-1)/3$, we shall
look at the variance of the number of edges and show that it is small
compared with $n(n-1)/3$. The way this will be done is by writing the
number of edges as the sum of indicator variables, one for each edge.
Evaluating the variance then involves considering the various possible values
of the expectation of $X_{ij}X_{k\ell}$ for various possible $i,j,k,\ell$.
We shall see that in most cases they are independent, so the total variance
is really rather small.
\begin{Theorem}
Almost all graphs in $\Delta_{n}$ have $n^{2}/3+ o(n^{2})$ edges.
More formally, letting $X$ denote the number of edges in the graph,
$$\lim_{n\rightarrow\infty}P\left(X=\frac{n^{2}}{3}(1+o(1))\right)=1.$$
\end{Theorem}
{\bf Proof.} For every pair $(i,j)$ with $1\leq i,j\leq n$ and $i\neq j$,
let the indicator random variable $X_{ij}$ be equal to $1$ if $i\sim j$ and
be equal to $0$ otherwise. (In other words, $X_{ij}$ takes the value of $1$,
if the intervals corresponding to the vertices $i$ and $j$ intersect, so
that the edge is present).

Now let $X =\Sigma_{1\leq i<j\leq n} X_{ij}$. Then each edge which is present
contributes 1 to this sum, and those absent contribute 0, thus $X$ counts
the total number of edges in the graph.

As we previously saw, $P(X_{ij}=1)=2/3$. Because $X_{ij}$ is an indicator
variable (i.e. takes only the values 1 and 0) we have
\begin{eqnarray*}
{\bf E}(X_{ij})=\sum_{r=0}^{\infty}rP(X_{ij}=r)=0+1\times P(X_{ij}=1)
=P(X_{ij}=1)
\end{eqnarray*}
Thus
\begin{eqnarray*}
{\bf E}(X)={\bf E}(\sum_{1\leq i<j\leq n}X_{ij})\\
=\sum_{1\leq i<j\leq n}{\bf E}(X_{ij})\\
=\sum_{1\leq i<j\leq n}P(X_{ij}=1)\\
=\frac{n(n-1)}{2}\frac{2}{3}\\
=\frac{n(n-1)}{3}
\end{eqnarray*}
To estimate the number more precisely, we need also to have some grip on
the variance of $X$, namely
\begin{eqnarray*}
\mbox{Var}(X)={\bf E}(X^{2})-({\bf E}(X))^{2}\\
={\bf E}(\sum_{1\leq i<j\leq n}X_{ij}\sum_{1\leq k<\ell\leq n}X_{kl})
-(\frac{n(n-1)}{3})^{2}
\end{eqnarray*}
(the last term is because we just worked out ${\bf E}(X)$)! Now we consider
various cases for $i,j,k$ and $\ell$, counting - or at least estimating -
carefully how many such cases there are.
\\
\\
{\bf Case 1.} $i,j,k,\ell$ are all distinct. As $i<j$ and
$k<\ell$, there are ${n\choose 2}{(n-2)\choose 2}/2$ possibilities (as the
first two are chosen from all $n$, the second two from the remaining
$(n-2)$. Thus we get $n(n-1)(n-2)(n-3)/4$ cases. In each of them,
${\bf E}(X_{ij}X_{kl})={\bf E}(X_{ij}){\bf E}(X_{kl})= 2/3\times 2/3=4/9$.
This is because the two edges are independent, as they have no vertex in
common.
\\
\\
{\bf Case 2.} There is some overlap between set $\{i,j\}$ and $\{k,\ell\}$.
There are however only a small number of these. Indeed since at least
two of $i,j,k$ and $\ell$ have to be equal, we are only choosing (at most)
three numbers. Thus there are at most ${n\choose 3}$ such cases. In each of
these cases, the contribution to the sum made by ${\bf E}(X_{ij}X_{kl})$
is at most 1, simply because $X_{ij}$ and $X_{kl}$ are $\leq 1$.

Then, calculating, we get
\begin{eqnarray*}
\mbox{Var}(X)={\bf E}(\sum_{1\leq i<j\leq n}X_{ij}\sum_{1\leq k<\ell\leq n}
X_{kl})-(\frac{n(n-1)}{3})^{2}\\
\leq \frac{n(n-1)(n-2)(n-3)}{4}\frac{4}{9}\mbox{~from~Case~1}\\
+ {n\choose 3}\times 1\mbox{~from~Case~2}\\
-(\frac{n(n-1)}{3})^{2}\\
=\frac{n(n-1)(n-2)(n-3)}{9}+{n\choose 3}-\frac{n^{2}(n-1)^{2}}{9}\\
=\frac{n(n-1)}{9}[(n-2)(n-3)+\frac{9(n-2)}{6}-n(n-1)]
\end{eqnarray*}
The above is, since the two terms in the square bracket involving $n^{2}$
cancel with each other, of the form
\begin{eqnarray*}
\frac{n(n-1)}{9}[Cn+D]\mbox{~for~suitable~constants~}C,D\\
\leq En^{3}
\end{eqnarray*}
for a suitable constant $E$.

Now we use Chebyshev's inequality, which says that, for any random
variable $X$,
$$P(\vert X-{\bf E}(X)\vert \geq \epsilon \leq
\frac{\mbox{Var}(X)}{\epsilon^{2}}.$$ See \cite{math}. Here we thus
deduce that, for any $\epsilon=cn^{2}$, we have
\begin{eqnarray*}
P(\vert X-{\bf E}(X)\vert \geq cn^{2} \leq
\frac{\mbox{Var}(X)}{(cn^{2})^{2}}\\
\leq \frac{En^{3}}{c^{2}n^{4}}
\end{eqnarray*}
which, for large $n$, tends to 0. Thus with probability tending to 1,
we do indeed get $(1+o(1))n^{2}/3$ edges. $\bullet$

It is perhaps worth noting that it is natural to compare a random
interval graph with another well-known model of random graphs. This
is the so-called Erd\H{o}s-R\'{e}nyi model $G(n,p)$ where there are
$n$ vertices and each edge arises with probability $p$ {\it
independently} of all other edges. This model is discussed in great
detail in \cite{bollobas}. Clearly the most reasonable such model to
compare random interval graphs with is $G(n,2/3)$ as we know that in
random interval graphs $2/3$ is the probability of each edge
arising. The above result does have a simple analogue for
Erd\H{o}s-R\'{e}nyi random graphs, namely that an Erd\H{o}s-R\'{e}nyi
random graph has about $n^{2}/3$ edges. The proof in this case is
much simpler: Indeed the law of large numbers, \cite{wiki} says that
the number of edges, divided by the total number of possible edges
$n(n-1)/2$, is close to the expectation of any one of the
indicators, namely $2/3$. Thus the number of edges is close to
$n(n-1)/3$. The reason why this case is so much easier is that the
edges in the Erd\H{o}s-R\'{e}nyi graph are independent, so that
standard results like the law of large numbers, \cite{wiki} apply to
them. (It will emerge later that the edges are not independent, when
we show that various things in a random interval graph usually take
very different values from their values in $G(n,2/3)$.)
\pagebreak

\noindent
\section{Degrees in Random Interval Graphs}
\subsection{Introduction}
We now turn our attention to the more detailed distribution of where the
edges
are. We shall consider the degrees of vertices and show that degrees are
much more spread out that in the Erd\H{o}s-R\'{e}nyi model $G(n,2/3)$.
Indeed in the Erd\H{o}s-R\'{e}nyi model $G(n,2/3)$, we shall see that, for all
$\epsilon>0$,
$$\lim_{n\rightarrow\infty}P\left(\mbox{all~vertices~have~degree~between~}
\frac{2n(1-\epsilon)}{3}\mbox{~and~}\frac{2n(1+\epsilon)}{3}\right)=1 .$$
That is, almost all degrees are about $2n/3$. However in random interval
graphs, we shall see in various ways that degrees are much more
\lq spread out\rq . For example, the probability that there is a vertex
of degree $n-1$ (the maximum possible degree) is $2/3$.

\subsection{Degrees of Graphs and some results}
We recall first a definition from \cite{iliop}.
\begin{Definition}
The degree of a node of a graph is the
number of vertices which are adjacent to this vertex. If
$v$ is a vertex, then the degree of $v$ is denoted by
\emph{deg(v)}.

For a graph $G$, we define $\Delta(G)$ to
be the maximum degree: that is,
$$\Delta(G)=\max_{1\leq i\leq n}d(v_{i}).$$
Similarly we define the minimum degree
$$\delta(G)=\min_{1\leq i\leq n}d(v_{i}).$$
\end{Definition}
We saw in the last section that a random interval graph has about
$n^{2}/3$ edges. Our next task is to show how these edges are spread out.
We will see a theorem from \cite{Schein} about
the degrees of random interval graphs.

\begin{Theorem}
Let $G\in\Delta_{n}$ and $v\in V(G)$. For a known $x\in [0, 1]$, we
have for $x\geq 1/2$
$$\lim_{n\rightarrow\infty}P(d(v)\leq xn) = 1-(1-x)\frac{\pi}{2}$$
and for $x<1/2$ we have
$$\lim_{n\rightarrow\infty}P(d(v)\leq xn) = 1-(1-x)(\pi /2 -
2\cos^{-1}[1/\sqrt{2-2x}])-\sqrt{1-2x}$$
\end{Theorem}
Note what the theorem means. It says that, for example, taking $x$ to be
(say) 0.01, the probability that there is a vertex of degree $\leq 0.01n$
is, in the limit as $n\rightarrow\infty$, strictly positive. This is
very different from what happens in $G(n,2/3)$ where, as mentioned
earlier, all the degrees are close to $2n/3$. Similarly it says that the
probability that a vertex does not have degree $\leq 0.99n$ - that is, that its
degree is at least $0.99n$ - is, in the limit, a
non-zero number. So the degrees are indeed much more spread out than
in $G(n,2/3)$.

The proof will rely on the following geometrical lemma. We give a rather
detailed proof of the Lemma
as no details are provided in Scheinerman's article \cite{Schein}. Some
of the details here were suggested to me by Dr. Penman \cite{penman}.
\begin{Lemma}
For an interval $I=[x,z]\subset [0,1]$ let
the radius of $I$, $\rho (I)$ be $\sqrt{a^{2}+(1-b)^{2}}$,
where $a=\min\{x,z\}$ and $b=\max\{x,z\}$. Assume that $x$ and $z$ are
independent, and uniformly distributed on $[0,1]$. Then, for $y\leq 1/2$
$$P(\rho^{2}(I)\leq y)= \frac{y\pi}{2}.$$
When $y>1/2$, we have that
$$P(\rho^{2}(I)\leq y) = y(\frac{\pi}{2}-2cos^{-1}[\frac{1}{\sqrt{2y}}])
+\sqrt{2y-1}.$$
\end{Lemma}

{\bf Proof.} The required probability is the probability that
$\rho(I)\leq \sqrt{y}$, which is the size of the set of points in
the square $[0,1]^{2}$ which are within a distance $\sqrt{y}$ from
(at least) one of the points $(0,1)$ or $(1,0)$. Let us consider all
the possible cases for $0\leq y\leq 1$. The equation of the circle
centered at $(1,0)$ is $({x-1}^{2})+z^{2}=y$. In the second circle,
centered at point $(0,1)$ we have the equation $x^{2}+({z-1}^{2})=y$.
So if they intersect at $(x,z)$, then
\begin{eqnarray*}
(x-1)^{2}+z^{2}=x^{2}+(z-1)^{2}=y\\
\Rightarrow
x^{2}-2x+1+z^{2}=x^{2}+z^{2}-2z+1\\
\Rightarrow -2z=-2x\Rightarrow z=x\\
\end{eqnarray*}
Thus, these points satisfy:
\begin{eqnarray*}
(x-1)^{2}+x^{2}=y\\
\Rightarrow 2x^{2}-2x+1=y\\
\Rightarrow x^{2}-x+\frac{1-y}{2}=0\\
\end{eqnarray*}
For such points to exist, we need that the discriminant is greater
than $0$. Thus,
\begin{eqnarray*}
\Rightarrow (-1)^{2}-4\frac{(1-y)}{2}\geq 0\\
\Rightarrow 1-2(1-y)\geq 0\\
\Rightarrow (2y-1)\geq 0\\
\Rightarrow y\geq\frac{1}{2}\\
\end{eqnarray*}
Now let consider the two cases for the value of the radius.
\\
\\
{\bf Case 1.} Let $y<\frac{1}{2}$, then $\sqrt{y}<1/\sqrt{2}$ is the
radius of the two circles centered at either the point $(0,1)$ or
$(1,0)$. Because $y<\frac{1}{2}$, these two circles do not intersect
by the above analysis. Then the area we want is two quarters of
disjoint circles which have radius equal to $\sqrt{y}$. Having in
mind that the area of a circle with radius $\chi$ is equal to
$\pi\times\chi^{2}$, then the area of a quarter of circle is
$\frac{1}{4}\pi\times\sqrt{y}^{2}=\frac{1}{4}\pi\times y$. So the
required area is equal to $\frac{2\pi\times y}{4}=\frac{\pi\times
y}{2}$, as stated.
\\
\\
{\bf Case 2.} $y\geq1/2$. In this case, the circles have
intersection and by the above we have for $x=z$, that:
\begin{eqnarray} x^{2}-x+\frac{(1-y)}{2}=0\\
\Rightarrow (x-\frac{1}{2})^{2}+\frac{(1-y)}{2}-\frac{1}{4}=0\\
\Rightarrow (x-\frac{1}{2})^{2}+\frac{(1-2y)}{4}=0\\
\Rightarrow (x-\frac{1}{2})^{2}+(\frac{1}{4}-\frac{y}{2})=0\\
\Rightarrow x=\frac{1}{2}\pm \sqrt{\frac{y}{2}-\frac{1}{4}}
\end{eqnarray}
In this case, as previously noted, the two circles intersect in two
points, $A$ and $B$ say. We consider the triangle with vertices
$(0,1)$, $A$ and $B$ and the triangle with vertices $(1,0)$, $A$ and
$B$, also the arcs of the two circles between $A$ and $B$. Let
$\theta$ be the angle formed at $[0,1]$ and $[1,0]$ by the two
triangles (by symmetry, the two angles are the same). The area of
each sector of each circle is equal to
$$\pi\times y\frac{\theta}{2\pi}=\frac{\theta y}{2}.$$
Also the area of the triangle is $\frac{1}{2}\times|CA|\times|CB|\times
\sin (\theta)$, where $C$ Is the point $(0,1)$ and A,B are the
intersection points, between the two circles.

Then the total area of intersection of the two sectors is:
\begin{eqnarray*}
2[\frac{\theta\times y}{2}-\frac{1}{2}|CA||CB|\sin\theta]\\
=\theta y-|CA||CB|\sin\theta
\end{eqnarray*}

To obtain the angle $\theta$, recall that the
vertices are: $C=(0,1)$,
$B=(\frac{1}{2}+\sqrt{\frac{y}{2}-\frac{1}{4}},
\frac{1}{2}+\sqrt{\frac{y}{2}-\frac{1}{4}})$
and
$A=(\frac{1}{2}-\sqrt{\frac{y}{2}-\frac{1}{4}},\frac{1}{2}-\sqrt{\frac{y}{2}
-\frac{1}{4}})$.

Now $|AB|^{2}=|AC|^{2}+|CB|^{2}-2|AC||CB|cos\theta$ by the cosine
rule, \cite{trig}. The vectors are
\begin{eqnarray*}
\overrightarrow{AB}=(2\sqrt{\frac{y}{2}-\frac{1}{4}},
2\sqrt{\frac{y}{2}-\frac{1}{4}})\\
\overrightarrow{CB}=(\sqrt{\frac{y}{2}-\frac{1}{4}}+\frac{1}{2},
\sqrt{\frac{y}{2}-\frac{1}{4}}-\frac{1}{2})\\
\overrightarrow{AC}=(-\sqrt{\frac{y}{2}-\frac{1}{4}}+\frac{1}{2},
-\sqrt{\frac{y}{2}-\frac{1}{4}}-\frac{1}{2})\\
\Rightarrow ÝABÝ^2=[2\sqrt{\frac{y}{2}-\frac{1}{4}}]^{2}
+[2\sqrt{\frac{y}{2}-\frac{1}{4}}]^{2}\\
=8(\frac{y}{2}-\frac{1}{4})=4y-2
\end{eqnarray*}
Similarly,
\begin{eqnarray*}
\vert CB \vert^{2}\\
=(\sqrt{\frac{y}{2}-\frac{1}{4}}+\frac{1}{2})^{2}
+(\sqrt{\frac{y}{2}-\frac{1}{4}}-\frac{1}{2})^{2}\\
=2[\frac{y}{2}-\frac{1}{4}]+\frac{1}{2}=y
\end{eqnarray*}
and then $|CA|^{2}=y$ just by symmetry. Then the cosine rule,
\cite{trig} becomes
\begin{eqnarray*}
(4y-2)=y+y-2\sqrt{y}\sqrt{y}\cos\theta\\
\Rightarrow 4y-2=2y-2y\cos\theta\\
\Rightarrow \cos\theta=\frac{2y-4y+2}{2y}\\
\Rightarrow \theta=\cos^{-1}[-1+\frac{1}{y}]\\
\end{eqnarray*}
(Note that, as $y\geq 1/2$, $-1+1/y\leq 1$ as required). Then, the
formula for the area of the intersection of the two quarter-circles is
(note that here we use $\cos^{-1}$ to mean the inverse function to
$\cos$, what many people call $\arccos$: in particular, $\cos^{-1}(x)$ does
{\it not} mean $1/\cos(x)$)
\begin{eqnarray*}
y\cos^{-1}[-1+\frac{1}{y}]-\sqrt{y}\sqrt{y}\sin[\cos^{-1}(-1+\frac{1}{y})]\\
=y\cos^{-1}[-1+\frac{1}{y}]-y\sqrt{1-(-1+\frac{1}{y})^{2}}\\
=y\cos^{-1}[-1+\frac{1}{y}]-\sqrt{y^{2}-(-y+1)^{2}}\\
=y\cos^{-1}[-1+\frac{1}{y}]-\sqrt{2y-1}
\end{eqnarray*}
To explain the working in the last paragraph; We used
$\cos^{2}\theta+\sin^{2}\theta=1$, \cite{trig} and also the
fact that in our situation the angle $\theta$ is clearly between 0
and $\pi/2$ so that both $\cos$ and $\sin$ are positive, with
the result that $\sin(x)=\sqrt{1-\cos^{2}(x)}$. Thus
$\sin[\cos^{-1}(-1+\frac{1}{y})]$ is equal to
$\sqrt{1-\cos^{2}(\cos^{-1}(-1+\frac{1}{y}))}$ which is of course
equal to $\sqrt{1-(-1+\frac{1}{y})^{2}}$.

What we have just worked out is the formula for the area of the
intersection of the two quarter-circles. Thus the shaded area in
Scheinerman's picture is the areas of the two individual quarter
circles minus the area of their intersection, which of course is
\begin{eqnarray*}
\frac{\pi y}{2}-\left(y\cos^{-1}[-1+\frac{1}{y}]-\sqrt{2y-1}\right)
=\frac{\pi y}{2}-y\cos^{-1}[-1+\frac{1}{y}]+\sqrt{2y-1}
\end{eqnarray*}
It only remains to confirm that this formula we have just derived is the
same as the one given in Scheinerman's article, namely
$$y(\frac{\pi}{2}-2cos^{-1}[\frac{1}{\sqrt{2y}}])
+\sqrt{2y-1}.$$
For this in turn it suffices to show that
$$2\cos^{-1}[\frac{1}{\sqrt{2y}}]=\cos^{-1}[-1+\frac{1}{y}].$$
To see this, take cosines of both sides, and recalling the identity
$\cos(2x)=2\cos^{2}(x)-1$, \cite{trig} we see the left-hand side is:
\begin{eqnarray*}
2\cos^{2}(\cos^{-1}[\frac{1}{\sqrt{2y}}])-1\\
=2[\frac{1}{\sqrt{2y}}]^{2}-1\\
=-1+\frac{1}{y}
\end{eqnarray*}
as required. $\bullet$
\\
\\
{\bf Proof of Scheinerman's theorem.}  (see Theorem 4.2 in
\cite{Schein}: we provide some more details). Let $v$ corresponds to
vertex $1$: no generality is lost by this, as no vertex is favored
by the set-up. Let $J_{i}$ be the interval assigned to vertex $i$,
for each $i =2\ldots n$. Then let $I_{i}=1$, if $1$ is adjacent to
$i$ and be equal to $0$ otherwise. Thus $X= \sum_{i=2}^{n} I_{i}$ is
the degree of the vertex $1$.

Suppose now that $\rho(I_{1})=r$ is fixed. We claim that then:
$$p=P(I_{i}=1\vert \rho(I_{1})=r)=1-r^{2}.$$
When we have proved this, it is then obvious that, conditional on
$\rho(I_{1})=r$, the expectation of $X$ is $(n-1)p$ and its variance
(again conditioned on the value of $\rho(I_{1})$) is $(n-1)p(1-p)$
since the $I_{i}$ are independent of each other given $I_{1}$ and
the radius. Indeed, if $i\neq j$, then the two random variables
$X_{i}$ and $Y_{i}$ giving the two endpoints of $I_{i}$ are
independent of $X_{j}$ and $Y_{j}$ giving the endpoints of $I_{j}$.

Thus we can use Chebyshev's inequality again on the random variable
$X$, \cite{math} and deduce that:
\begin{eqnarray*}
P(|X-(n-1)p|\geq n^{2/3})\leq \frac{\mbox{Var}(X)}{n^{4/3}}\\
\leq \frac{np(1-p)}{n^{4/3}}\rightarrow 0\\
\end{eqnarray*}
Then, $d(1)=np+o(n)$, under the hypothesis that $\rho(I_{1})=r$.
Thus for any $\epsilon>0$ we have:
\begin{eqnarray*}
P(d(1)\leq xn)=\left \{
\begin{array}{r}
1-o(1)\mbox{~for~}r<\sqrt{1-x}-\epsilon\\
o(1)\mbox{~for~}r>\sqrt{1-x}+\epsilon
\end{array}
\right .
\end{eqnarray*}
using our formula for the value of $p$ for a given value of $r$. So now
we need to remove the conditioning on the value of $r$, which we do in
the usual manner:
\begin{eqnarray*}
P(d(1)\leq xn)=\int_{0}^{1}P(d(1)\leq xn \vert p(I_{1})=r)
dP(\rho(I_{1})\leq r)\mbox{~by~the~law~of~total~probability}
\end{eqnarray*}
But we have just worked out the distribution function for the
probability that $\rho(I_{1})\leq \sqrt{y}$ in the Lemma. So this
is:
\begin{eqnarray*}
[1-o(1)]P(\rho^{2}(I_{1})\leq 1-x)+\epsilon O(1)+o(1)\\
\Rightarrow P(\rho^{2}(I_{1})\leq 1-x).
\end{eqnarray*}
Now the result follows using the Lemma, with $y$ replaced by $1-x$.

Thus the only thing that remains to be proved is that, if
$\rho(I_{1})=r$, then the probability that any other interval
intersects is $1-r^{2}$. To obtain this, recall that in
\cite{Schein}, the radius of an interval $[x,y]\subset [0,1]$ is
equal to $\sqrt{a^{2}+(1-b)^{2}}=r$, where $a=\min\{x,y\}$ and
$b=\max\{x,y\}$. Also, we should remind ourselves that the $(2n)$
random points follow the uniform distribution in $(0,1)$. Let
$I_{1}=[a,b]$. We calculate the possibility of the existence of $i$
interval, denoted by $I_{i}=[X_{i},Y_{i}]$, which does not intersect
with $I_{1}$. To happen this, both $X_{i}$, $Y_{i}$ must be smaller
than $a$ or both must be greater than $b$. The first possibility is:
$P(X_{i}, Y_{i}\leq a)=\frac{a-0}{1-0}\frac{a-0}{1-0}=a^{2}$, since
they are independent, uniformly distributed random variables in
$(0,1)$. The second possibility is equal to:
\begin{eqnarray*}
P(X_{i}, Y_{i}\geq b)\\
=P(X_{i}\geq b)P(Y_{i}\geq b)\\
=(1-P(X_{i}\leq b))(1-P(Y_{i}\leq b))\\
=(1-\frac{b-0}{1-0})(1-\frac{b-0}{1-0})\\
=(1-b)(1-b)\\
=(1-b)^{2}.
\end{eqnarray*}
Thus, the probability of $I_{1}\cap I_{i}=\emptyset$, for $i=2,
3\ldots n$ is equal to $a^{2}+(1-b)^{2}=r^{2}$. Hence the complement
probability ($I_{1}$ and $I_{i}$ have a non-null intersection) is
$1-r^{2}$.  $\bullet$

\subsection{Maximum degree}
In the last section, we saw there is a non-zero probability that a
vertex has very high degree (e.g. at least $0.99n$). In this section, we
first sharpen this. We then give a result from
\cite{JSTOR} which \cite{Schein} attempted
to prove, but did not succeed. That result is, that with probability
$2/3$, there is a vertex in a random interval graph whose degree is
$n-1$. This is of course the largest possible degree any vertex in
the graph can have.

The first result (which is Theorem 4.4. in Scheinerman) is easy.
\begin{Theorem}
In a random interval graph, let $\omega_{n}$ be any function which tends
to infinity with $n$. (One should think of it as doing so very slowly).
Then
$$\lim_{n\rightarrow\infty}P\left(\mbox{a~random~interval~graph~has~}
\Delta\geq n-\omega(n)\right)=1.$$
\end{Theorem}
{\bf Proof.} Some details of that proof were suggested to me by Dr.
Penman \cite{penman}, as Scheinerman's proof for that Theorem is
rather short. Let $x=\frac{1}{2}\sqrt{\frac{\omega_{n}}{n}}$. Then
\begin{eqnarray*}
P(\geq n-\frac{\omega_{n}}{2}+o(\omega_{n}) \mbox{~intervals~intersect~}[x,1-x])\\
=1-P(\geq n-\frac{\omega_{n}}{2}+o(\omega_{n})\mbox{~intervals~don't~intersect}[x,1-x])\\
=1-P(\mbox{interval~}I_{1}\mbox{~doesn't~intersect~}[x,1-x])^{n}\\
=1-[(\frac{1}{2}
\sqrt{\frac{\omega_{n}}{n}}+\frac{1}{2}\sqrt{\frac{\omega_{n}}{n}})^{2}]^{n}
\end{eqnarray*}
as the probability a random interval does not intersect $[x,1-x]$ is the
probability that both ends are less than $x$ (probability $x\times x=x^{2}$,
using the independence of the ends) or that both ends are greater than $1-x$
(probability $[1-(1-x)]^{2}=x^{2}$). The above is
\begin{eqnarray*}
=1-(\frac{\omega_{n}}{n})^{n}\\
\rightarrow 1\mbox{~as~$n\geq \omega_{n}$~when~$n\rightarrow\infty$}
\end{eqnarray*}
Thus at least $n-\omega_{n}/2+o(\omega_{n})$ intervals intersect $[x,1-x]$.
Our proof will now be complete if we can show that, with probability tending
to 1 as $n\rightarrow\infty$, there is at least one of the random intervals
which contains $[x,1-x]$, as then such an interval will be a vertex of degree
$\geq n-\omega(n)$.

To this end, recall the fact that random intervals are independent, as all
the possible orderings of their endpoints are equally likely. Let $X_{i}$ be
an indicator variable, denoting whether the interval $I_{i}=[A_{i},B_{i}]$
contains $[x,1-x]$. Then $X=\sum_{i=1}^{n}X_{i}$ is the number of our $n$
random intervals which contain $[x,1-x]$ and our objective is to show that
$X>0$ with probability tending to 1. To this end we use Chebyshev's inequality,
\cite{math}
\begin{eqnarray*}
P(|X-{\bf E}(X)|\geq t)\leq \frac{Var(X)}{t^{2}}\\
\Rightarrow P(X=0)\leq P(\vert X-{\bf E}(X)\vert \geq {\bf E}(X))\\
\leq \frac{Var(X)}{({\bf E}(X))^{2}}
\end{eqnarray*}
Now $X_{i}$ is a Bernoulli variable and the various $X_{i}$ are independent
of each other, simply because the distinct intervals are independent of
each other. Thus $X$ is binomial, with parameters $n$ and the success probability
$a$. We thus have, by standard results about Binomial random variables,
\begin{eqnarray*}
{\bf E}(X)=na\mbox{~and~Var}(X)=\sum_{i=1}^{n}Var(X_{i})=na(1-a)\\
\end{eqnarray*}
and thus by Chebyshev, \cite{math}
\begin{eqnarray*}
P(X=0)\leq \frac{Var(X)}{({\bf E}(X))^{2}}\\
=\frac{na(1-a)}{n^{2}a^{2}}=\frac{(1-a)}{na}
\end{eqnarray*}
Thus if we can show that $P(X=0)$ tends to zero, this will show that
$P(X>0)$ tend to 1.

We need to calculate $a$. Recall
\begin{eqnarray*}
a=P(X_{i}=1)=P(\mbox{a~given random~interval~contains~}[x,1-x])\\
\end{eqnarray*}
Letting the two (random) ends of the interval be $A$ and $B$, we have
that this probability is
\begin{eqnarray*}
a=P(\{A\leq x\mbox{~and~}B\geq 1-x)\}\cup \{A\geq 1-x\mbox{~and~}B\leq x\})\\
=P(A\leq x\mbox{~and~}B\geq 1-x))+P(A\geq 1-x\mbox{~and~}B\leq x)
\end{eqnarray*}
as the two events involved are mutually exclusive. This in turn gives
\begin{eqnarray*}
a=P(A\leq x)P(B\geq 1-x)+P(A\geq 1-x)P(B\leq x\})\mbox{~by~independence}\\
=x(1-x)+x(1-x)=2x(1-x)\\
=2\frac{1}{2}\sqrt{\frac{\omega_{n}}{n}}
(1-\frac{1}{2}\sqrt{\frac{\omega_{n}}{n}})\\
=\sqrt{\frac{\omega_{n}}{n}}(1-\frac{1}{2}\sqrt{\frac{\omega_{n}}{n}})
\end{eqnarray*}
Hence, combining the above results,
\begin{eqnarray*}
P(X=0)\leq \frac{1-a}{na}\leq \frac{1}{na}\\
=\frac{1}{n\sqrt{\omega_{n}/n}\left(1-\sqrt{\omega_{n}/(2n)}\right)}\\
=\frac{1}{n^{1/2}\sqrt{\omega_{n}}(1-\sqrt{\omega_{n}/(2n)})}\\
\leq \frac{2}{n^{1/2}\sqrt{\omega(n)}}\mbox{~for~large~enough~n}
\end{eqnarray*}
using the fact that
$$1-\frac{1}{2}\sqrt{\frac{\omega_{n}}{n}}\geq \frac{1}{2}$$
for large enough $n$. Thus
\begin{eqnarray*}
P(X=0)\leq \frac{2}{n^{1/2}\omega(n)}\rightarrow 0\mbox{~as~}
n\rightarrow\infty
\end{eqnarray*}
and thus with probability tending to $1$, $X>0$, that is there is
some interval containing $[x,1-x]$ as required. $\bullet$
\\
\\
The final result is the following striking fact. We emphasize that,
unlike most of the results in this chapter, this has nothing to do
with a limit: it is an {\it exact} result, not depending on
the number of vertices. The proof comes from \cite{JSTOR}. Scheinerman
made a lot of effort in his paper \cite{Schein} to prove a result along
these lines, but did not quite succeed.
\begin{Theorem}
Let $G$ be a random interval graph. Then
$$P(\Delta(G)=n-1)=\frac{2}{3}.$$
\end{Theorem}
{\bf Proof.} Instead of calculating the moderately difficult
expression
\begin{eqnarray*}
1-4n(n-1)\int_{0}^{1}\int_{0}^{1-y}xy(1-x^{2}-y^{2}-2xy)^{n-2}-1dxdy\\
\end{eqnarray*}
writers \cite{JSTOR} used an efficient combinatorial proof, as we
shall see. They take random pairs of integers $1, 2,\ldots, 2n$.
Once the intervals are selected by some random pairing of the $2n$
numbers, they label the endpoints $A(1), B(1),\ldots, A(n-2),
B(n-2)$ in the following way. Let the endpoints $\{1,\ldots, n\}$ be
at the left side and respectively the endpoints $\{n+1,\ldots,
2n\}$, be at the right side. Let also $A(1)=n$ and $B(1)$ is its
mate. Suppose that we have assigned through $A(j)$, $B(j)$. We label
the next endpoints, by the following rules:
\\
\\
{\bf Case 1.} If $B(j)$ is on the left side, then let $A(j+1)$ be
the leftmost point on the right side that has not yet been labeled.
Let $B(j+1)$ be its mate.
\\
\\
{\bf Case 2.} If $B(j)$ is on the right side, then let $A(j+1)$ be
the rightmost point on the left side that has not yet been labeled.
Let $B(j+1)$ be its mate.
\\
\\
Endpoints are being labeled from the center outwards. Then,
if $A(j)<B(j)$, it is $A(j+1)<B(j)$\\
If $A(j)>B(j)$, it is $A(j+1)>B(j)$\\
Then, we will either have:
\begin{eqnarray*}
A(j)<A(j+1)\\
\mbox{~or~} A(j)>A(j+1)
\end{eqnarray*}

Furthermore, if $A(j)<B(j)$, then $A(j+1)<B(j+1)$\\
If $A(j)>B(j)$, then $A(j+1)>B(j+1)$\\

With this way, starting from the center labeling the endpoints, we
deduce that either an equal number of points have been assigned in
both sides, or two more points have been assigned on the left than
on the right side. Since the last endpoints assigned are $A(n-2)$
and $B(n-2)$, from the total number of points, which is equal to
$2n$, there are four remaining points unlabeled, namely $a<b<c<d$.
Having a specific ordering, it is considered all the possible ways
of pairing them to consist two random intervals. $a$ can be matched
with $b$, $c$ or $d$, with equal probabilities. Thus, we have $3$
possible cases. We easily observe that in two cases of pairing the
remaining points, the corresponding intervals intersect and in one
case they are disjoint. Let now $a$ and $b$ be on left and $c$ and
$d$ on the right. If $a$ is paired with $c$, then the random
interval $[a,c]$ meets all the others. Also the same happens when
$a$ is paired with $d$. This is because, we assumed the points $\{1,
2\ldots n\}$ lie on left and the remaining points $\{n+1\ldots 2n\}$
lie on the right and by construction of labeling them. On the other
hand, if $a$ is paired with $b$, then $[a,b]\cap[c,d]=\emptyset$.
Suppose that an interval $[e,f]$ intersects all the others. Also let
$A(j)=e$ and $B(j)=f$. In this case, where $a$ and $b$ are on the
left, $A(j)$ lies between $b$ and $c$. Thus, $[e,f]$ cannot
intersect both $[a,b]$ and $[c,d]$. Furthermore, consider the case
where only $a$ is on left. Since $[e,f]$ meets $[c,d]$ we have
$f>c$, hence $f=B(j)$. Again, by construction if $a$ is paired with
$b$, then $[a,b]\cap[c,d]=\emptyset$. On the other hand, in cases
where $a$ is paired with $c$ or $d$ the corresponding intervals
intersect. The probability of pairing $a$ with $c$ or $d$, in the
specific ordering of the four endpoints, which is $a<b<c<d$ is
$2/3$. Finally, we see that the probability in a family of $n$
random intervals, the maximum degree has value $n-1$ (i.e the
probability that an interval meets all the others) is $2/3$.
$\bullet$

\subsection{Minimum degree}
The previous subsection makes it clear that the maximum degree in a
random interval graph is much bigger than $2n/3$, which is roughly
the expected number of neighbors of each vertex. As we will see at
the end of this chapter, this is different from $G(n,2/3)$ where all
degrees are about $2n/3$. We now say a little about minimum degrees,
presenting a Theorem from \cite{Schein}. We omit the proof of this
result, as it is moderately difficult.
\begin{Theorem}
Let $k$ be a fixed, non-negative real number and $\delta$ denotes
the minimum degree of the graph. We have,
$$\lim_{n\rightarrow\infty}P(\delta<k\sqrt{n})=1-\exp\{-\frac{k^{2}}{2}\}$$
\end{Theorem}
Note that Theorem $4.1$ implies the result below, which is
noted in \cite{Schein}.
\begin{Corollary}
For every $\epsilon>0$ sufficiently small, almost all interval
graphs satisfy $\delta<\epsilon\cdot n$ and
$\Delta>(1-\epsilon)\cdot n$.
\end{Corollary}
This result is an immediate corollary of Theorem 4.1.
\subsection{Degrees of vertices in $G(n,2/3)$}
In the previous subsections, we studied the minimum and maximum
degrees of Random Interval Graphs. We now give, for contrast, the result
for the Erd\H{o}s-R\'{e}nyi model, where the probability of an edge arising
is constant, equal to $2/3$ and the edges are independent. Here it will
turn out that all the degrees are \lq about\rq\, $2n/3$.
\begin{Theorem}
In $G(n,\frac{2}{3})$ and for $\epsilon>0$ sufficiently \lq small\rq , all
vertices have degree between about $(\frac{2}{3}-\epsilon)n$ and
$(\frac{2}{3}+\epsilon)n$. More precisely,
\begin{eqnarray*}
\lim_{n\rightarrow \infty}P(\mbox{all~vertex~degrees}\in
[(\frac{2}{3}-\epsilon)n,(\frac{2}{3}+\epsilon)n])=1\\
\end{eqnarray*}
\end{Theorem}
Some details of this proof were suggested by \cite{penman}. Also,
for the proof of this Theorem, we use the \lq Large Deviations\rq\,  Lemma
from Scheinerman, \cite{Schein}.
\begin{Lemma}
If $p$ is constant and $\epsilon>0$, then
\begin{eqnarray*}
P(|X-np|\geq \epsilon
np)\leq\frac{a_{\epsilon}e^{-b_{\epsilon}pn}}{\sqrt{np}}
\end{eqnarray*}
Where $a_{\epsilon}$, $b_{\epsilon}$ are positive constants, which
depend only on $\epsilon$: that is, not on $n$ or $p$.
\end{Lemma}
{\bf Proof of Theorem 4.7.} Let the random variable $X$ denotes the
number of vertices having degrees not in the interval
$$[(\frac{2}{3}-\epsilon)n, (\frac{2}{3}+\epsilon)n].$$
$X$ is the sum
of $n$ independent random variables $X_{i}, i=1, 2\ldots n$ where
$X_{i}$ is 1 if vertex $i$ has degree not in the stated range and is zero
otherwise. We again aim to show that $X=\sum_{i=1}^{n}X_{i}$ is 0 with
probability tending to one as $n\rightarrow\infty$, and to prove this
we will use the fact that
$$P(X>0)=\sum_{i=1}^{\infty}P(X=i)\leq \sum_{i=1}^{\infty} iP(X=i)
\leq {\bf E}(X).$$
So we aim to show ${\bf E}(X)\rightarrow 0$, for which in turn it suffices,
as
\begin{eqnarray*}
{\bf E}(X)={\bf E}(\sum_{i=1}^{n}X_{i})=\sum_{i=1}^{n}{\bf E}(X_{i})\\
=n{\bf E}(X_{1})\mbox{~as~the~}X_{i}\mbox{~are~identically~distributed}\\
=nP(X_{1}=1)
\end{eqnarray*}
to show that $P(X_{1}=1)$ is $o(n)$.

But the degree of vertex $i$ in $G(n,\frac{2}{3})$ has binomial
distribution with parameters $n-1$ and $\frac{2}{3}$. Thus, using
the Large Deviations Result, \cite{Schein}
\begin{eqnarray*}
P[\mbox{degree~of~vertex~i}\notin[(\frac{2}{3}-\epsilon)n,(\frac{2}{3}+\epsilon)n]]\\
\leq \frac{a_{\epsilon}e^{-b_{\epsilon}\frac{2}{3}n}}{\sqrt{\frac{2n}{3}}}\\
\Rightarrow E(X_{i})=P(X_{i}=1)\times1+P(X_{i}=0)\times0=P(X_{i}=1)\\
E(X_{i})=P(X_{i}=1)\leq \frac{a_{\epsilon}e^{-b_{\epsilon}\frac{2}{3}n}}{\sqrt{\frac{2n}{3}}}\\
\Rightarrow E(X)
=\sum_{i=1}^{n}E(X_{i})\leq n\frac{a_{\epsilon}e^{-b_{\epsilon}\frac{2}{3}n}}
{\sqrt{\frac{2n}{3}}}\\
\end{eqnarray*}
And this is indeed $o(n)$ for large $n$, just because the
exponential terms converge more rapidly than polynomial. $\bullet$
\\
\\
{\bf Remark.} In fact the above argument can be sharpened quite substantially:
a random graph $G(n,2/3)$ has the property that there is an explicit constant
$C$ such that
\begin{eqnarray*}
\lim_{n\rightarrow\infty}P\left(\frac{2n}{3}-C\sqrt{n\log(n)}\leq \delta(G)
\leq \Delta(G)\leq \frac{2n}{3}+C\sqrt{n\log(n)}\right)=1.
\end{eqnarray*}
The proof is in \cite{bollobas} (but is much harder than the above proof).
\pagebreak

\noindent
\section{Cliques, independent sets and chromatic numbers in
Random Interval Graphs}
\subsection{Cliques}
In this section we study the clique number $\omega(G)$ of a random
interval graph $G$. Recalling from \cite{iliop}, the clique number
of a graph $G$, denoted as $\omega(G)$ is the order of the largest
complete subgraph of $G$. The basic result is from \cite{Schein},
Theorem 4.7.
\begin{Theorem}
The clique number of a random interval graph is usually about $n/2$.
More precisely,
$$\lim_{n\rightarrow\infty}P\left(\omega(G)=\frac{n}{2}+o(n)\right)=1.$$
\end{Theorem}
{\bf Proof.} First note that a maximum clique consists of a family
of intervals, each pair of which intersect. In our essay,
\cite{iliop} we showed that such a family of intervals has some
point $x$ say, which is in all the intervals. This result is known
as Helly's Theorem.

It is intuitively obvious that the point with the best chance of
being in several intervals is $x=1/2$: let us be more formal about
this now, using the argument from \cite{JSTOR}. The probability that
a random interval $[X,Y]$ does not contain $x$ is the probability
that both $X$ and $Y$ are less than $x$, which has probability
$x^{2}$, or the (exclusive) possibility that both are greater than
$x$, which is $(1-x)^{2}$: thus the probability that it does contain
$x$ is $1-x^{2}-(1-x)^{2}=2x-2x^{2}$. This is maximized for $x=1/2$,
as the derivative of $2x-2x^{2}$ is $2-4x$ which is zero exactly
when $x=1/2$, when the probability that an interval contains $x$ is
$2\cdot 1/2-2\cdot (1/2)^{2}=1/2$. In a different sense, since the
endpoints of the $n$ random intervals follow the the uniform
distribution in $[0,1]$, then with probability equal to $1$ will be
distinct. Moreover, the mean number of a variable, which follows the
uniform distribution is equal to $\frac{0+1}{2}=\frac{1}{2}$. As a
result, the expected number of intervals is indeed $n/2$. The number
of $n$ random intervals containing $1/2$ is binomial:
$\mbox{Bin}(n,1/2)$. That is because each random interval in
$(0,1)$, either it will contain $1/2$ with probability $p$ say, or
it will not contain it, with probability $1-p$. We have $n$
independent Bernoulli trials so the number of intervals containing
$1/2$ is indeed a Binomial random variable. This takes a value very
close to $n/2$, as we saw above: in particular, it is $n/2+o(n)$.
This gives a lower bound on the clique number.

We now need to show that it is not more than $n/2$.
Consider intervals of the form $I_{i}=[i/n^{2},(i+1)/n^{2}]$. We need
the following technical result (which Scheinerman calls the \lq
medium deviations lemma\rq\, for binomial random variables):
\\
\\
\begin{Lemma} If $X=\sum_{i=1}^{n}X_{i}$ where $X_{i}=1$ with
probability $p$ and is zero otherwise, and the $X_{i}$ are
independent, then for
$$1\leq h < \min \left\{\frac{np(1-p)}{10},\frac{(pn)^{2/3}}{2}\right\}$$
then
$$P\left(\vert X-np\vert\geq h\right)\leq \frac{\sqrt{np(1-p)}}{h}
\exp(-\frac{h^{2}}{2np(1-p)}).$$
\end{Lemma}
The above lemma comes from \cite{Schein}.
\\
\\
How does this help us? Let us work out the probability that more
than $n/2+o(n)$ of the $n$ random intervals intersect $I_{i}$. (We are
aiming to show that this does not happen). The
probability that a random interval $[X,Y]$ does not intersect
$I_{i}$ is the probability that both $X$ and $Y$ are less than
$i/n^{2}$ (which is $i^{2}/n^{4}$), or that both are bigger than
$(i+1)/n^{2}$ which has probability $(1-\frac{i+1}{n^{2}})^{2}$.
Thus the success probability (i.e the probability that it does
intersect $I_{i}$) is
\begin{eqnarray*}
1-\frac{i^{2}}{n^{4}}-\left(1-\frac{i+1}{n^{2}}\right)^{2}\\
=1-\frac{i^{2}}{n^{4}}-(1-\frac{2(i+1)}{n^{2}}+\frac{(i+1)^{2}}{n^{4}})\\
=\frac{2(i+1)}{n^{2}}-\frac{i^{2}}{n^{4}}-\frac{(i+1)^{2}}{n^{4}}
\end{eqnarray*}
Hence the number $X$ of the $n$ random intervals which intersect
$I_{i}$ is binomial with $n$ trials and success probability
$$\frac{2(i+1)}{n^{2}}-\frac{i^{2}}{n^{4}}-\frac{(i+1)^{2}}{n^{4}}.$$
We need to work out which value of $i$ maximizes this expression.
Let us treat the more general problem of which (continuous) value of
$x$ maximizes
$$f(x)=\frac{2(x+1)}{n^{2}}-\frac{x^{2}}{n^{4}}-\frac{(x+1)^{2}}{n^{4}}.$$
Differentiating, we get
\begin{eqnarray*}
f'(x)=\frac{2}{n^{2}}-\frac{2x}{n^{4}}-\frac{2(x+1)}{n^{4}}\\
=\left(\frac{2}{n^{2}}-\frac{2}{n^{4}}\right)-\frac{4x}{n^{4}}
\end{eqnarray*}
so the turning point is at $f'(x)=0$,
$$x=\frac{n^{4}}{4}\left(\frac{2}{n^{2}}-\frac{2}{n^{4}}\right)
=\frac{n^{2}-1}{2}.$$
Since $f"(x)<0$, this turning point is a maximum. Note that at it
we have
\begin{eqnarray*}
f(x)=\frac{n^{2}+1}{n^{2}}-\frac{(n^{2}-1)^{2}}{4n^{4}}
-\frac{(n^{2}+1)^{2}}{4n^{4}}\\
=1+\frac{1}{n^{2}}-\frac{1}{4}+\frac{1}{2n^{2}}-\frac{1}{4n^{4}}
-\frac{1}{4}-\frac{1}{2n^{2}}-\frac{1}{4n^{4}}\\
=\frac{1}{2}+\frac{1}{n^{2}}-\frac{1}{2n^{4}}\\
=\frac{1}{2}+\frac{2n^{2}-1}{2n^{4}}
\end{eqnarray*}
The point about this number is that it is very close to $1/2$.

We now (rather arbitrarily) consider 0.6. Note that
$n^{0.6}<(np)^{2/3}/2$ and $n^{0.6}<np(1-p)/10$ for large enough
values of $n$. The same holds with 0.6 replaced by 0.61. (This is
checked so that we can apply the medium deviations lemma in a minute).
We have
\begin{eqnarray*}
P(X\geq n/2+n^{0.6})\\
\leq P(\mbox{Bin}(n,1/2+(2n^{2}-1)/n^{4})\geq n/2+n^{0.6})\\
\end{eqnarray*}
since the success probability for the latter binomial is the largest
possible value of the success probability under the constraints.

This we now apply the medium deviations result to, taking $h=n^{0.6}$ and
$t=h/\sqrt{np(1-p)}$ where $p=1/2+(2n^{2}-1)/n^{4}$. Note that given
$\epsilon >0$, $p<1/2+\epsilon$ for all large enough $n$. Thus $p(1-p)$ is
$\geq (1/2+\epsilon)(1/2-\epsilon)=1/4-\epsilon^{2}$, and on the other
hand is $\leq 1/4$. Thus $t$ is between $\sqrt{1/4-\epsilon^{2}}n^{0.1}$
and $1/2 n^{0.1}$. Thus
\begin{eqnarray*}
\leq P(\mbox{Bin}(n,1/2)\geq n/2+n^{0.6})\leq \frac{1}{t}\exp(\frac{-t^{2}}{2})
\end{eqnarray*}
and this is, by the above estimates for $t$,
\begin{eqnarray*}
\leq \frac{1}{\sqrt{1/4-\epsilon^{2}}n^{0.1}} e^{-(1/4-\epsilon^{2})n^{0.2}/2}
\rightarrow 0\mbox{~as~}n\rightarrow\infty
\end{eqnarray*}
and this completes the proof. $\bullet$
\\
\\
The following simple consequence is also contained in Theorem 4.7 of
\cite{Schein}.
\begin{Corollary}
The chromatic number of a random interval graph satisfies
$$\lim_{n\rightarrow\infty}P\left(\chi(G)=\frac{n}{2}+o(n)\right)=1.$$
Informally: it is usually approximately $n/2$.
\end{Corollary}

{\bf Proof.} Interval graphs are perfect, and all perfect graphs
have $\omega(G)=\chi(G)$, as we saw in our essay \cite{iliop}. The
result now follows from the previous theorem. $\bullet$
\\
\\
Here is a comparison of the result with what happens for $G(n,2/3)$,
though we do not give the proof. It turns out that the clique number
of $G(n,2/3)$ is (with probability tending to 1 as
$n\rightarrow\infty$) about $2\log_{3/2}(n)$, which of course is far
smaller than for the random interval graphs. Also, in
Erd\H{o}s-R\'{e}nyi, the chromatic number is (again, with
probability tending to 1 as $n\rightarrow\infty$) about
$n/(2\log_{3}(n))$, which is of course much larger than the clique
number. In particular, the Erd\H{o}s-R\'{e}nyi graph is very far
from being a perfect graph, since the chromatic number is (for large
$n$) so much larger than the clique number. Note that in both
the case of the clique number and the case of hte chromatic number,
we get noticeably larger answers for random interval graphs than we do
for Erd\H{o}s-R\'{e}nyi graphs. We refer to \cite{bollobas} for detailed
statements of these two results.

Here we present a reformulation from \cite{JSTOR} of the result
we have just proved about the clique number.
\begin{Theorem}
Let the random variable $A_{n}$ denotes the size of the largest set
of pairwise intersecting intervals in a family of $n$ random
intervals. There exists a function $f(n)$ such that:
\begin{eqnarray*}
\lim_{n\rightarrow\infty} \frac{f(n)}{n}=0\\
\mbox{~and~} \lim_{n\rightarrow\infty}P(\frac{n}{2}-f(n)\leq A_{n}\leq \frac{n}{2}+f(n))=1\\
\end{eqnarray*}
\end{Theorem}

\subsection{Independent sets in random interval graphs}
Recall from \cite{iliop} that an independent set in a random
interval graph is the same as a chain in the random interval order
associated with it. We thus investigate the length of the longest
chain in the partial order, so as to obtain a result on the size of
the largest independent set. Let a graph with $n$ vertices where the
set of edges is equal to the emptyset, then the corresponding random
interval graph consists of $n$ disjoint intervals. Thus, the size of
the maximum chain is equal to $n$. The ordered intervals form a
chain, in the sense that they have a null intersection. On the other
hand, in a complete graph where all the edges are existent, the
chain is equal to the null set, as all the intervals intersect.
Again, this material is based on \cite{JSTOR} sharpening results in
\cite{Schein}. Recalling from \cite{iliop} a partially ordered set
is a non-empty set, which has the mathematical property of order.
Taking randomly two elements, say $x$ and $y$, then if we have $x<y$
or $y<x$ we define these elements as being \emph{comparable} forming
a \emph{chain}. In other case, we define them as
\emph{incomparable}, forming an \emph{antichain}. In the case of a
random interval order, where for two intervals $[a,b]$ and $[c,d]$,
we say $[a,b]\prec [c,d]$ if $b<c$. Thus an independent set in the
random interval graph corresponds to a set of non-intersecting
intervals.
\begin{Theorem}
Let $Y_{n}$ denote the maximum number of pairwise disjoint intervals
in a family of $n$ random intervals. Then,
$$\lim_{n\rightarrow\infty}\frac{Y_{n}}{\sqrt{n}}=\frac{2}{\sqrt{\pi}}$$
in probability.
\end{Theorem}
The proof of this Theorem comes from \cite{JSTOR}.
\\
\\
{\bf Proof.} We generate a Poisson process with intensity 1
in the upper right quadrant. Thus, since the probability
function for a discrete Random Variable $X$, which follows the
Poisson distribution is $\frac{e^{-\lambda}\lambda^{x}}{x!}$, where
$\lambda$ is the expected value of $X$, the probability that,
for any positive number $s$, the region $\{(x,y):0\leq x, y\leq s\}$ does
not contain any point of the process is equal to $e^{-s^{2}}$.

We are going to choose an infinite chain of points of the process.
Let $C=\{(\ell_{1},u_{1}), (\ell_{2},u_{2}),\ldots \}$, the points
being chosen as follows. Let $(\ell_{1}, u_{1})$ be the point that
minimizes $\max\{\ell_{1}, u_{1}\}$ and thereafter
$(\ell_{k},u_{k})$ is the point above $(\ell_{k-1},u_{k-1})$ that
minimizes $\max\{\ell_{k},u_{k}\}$. Thus the points of the chain are
chosen subject to the restrictions that they are monotonically
increasing and the difference between every two \lq neighboring \rq
points of the chain minimum. Thinking of each point in the chain as
defining an interval, it is easy to see that thus this chain is
built up from the bottom by always choosing the next interval to be
the one with least possible upper endpoint. It is not hard to check
by induction the intuitively reasonable claim that, in any finite
collection of intervals, this chain will have the maximum possible
length.

Then, if $S$ is a variable, whose value is $\max(x_{1},y_{1})$ the
mass density function of $S$ is given by:
\begin{eqnarray*}
f(s)=\frac{d}{ds}(1-e^{-s^{2}})=2se^{-s^{2}}.
\end{eqnarray*}
This is because we have
\begin{eqnarray*}
P(\max(x_{1},y_{1})\leq s)=1-P(x_{1}\geq s\mbox{~and~}y_{1}\geq s)\\
=1-P(\mbox{the~square~}(0,0),(0,s),(s,0)\mbox{~and~}(s,s)
\mbox{~contains~no~point~of~the~process})\\
=1-e^{-s^{2}}\mbox{~as~observed~above,~using~that~the~intensity~is~}1
\end{eqnarray*}
Thus $F(s)=P(S\leq s)=1-e^{-s^{2}}$. Thus its density $f(s)$ is the
derivative of this with respect to $s$, which is indeed as stated.

Therefore, we have
\begin{eqnarray*}
E(S)=\int_{0}^{\infty}2se^{-s^{2}}s.ds\\
=\int_{0}^{\infty}t^{\frac{1}{2}}e^{-t} dt\\
=\Gamma(3/2)=\frac{\sqrt{\pi}}{2}\\
\end{eqnarray*}
We now claim the differences
$$X_{1}=\max(\ell_{1},u_{1})-0, X_{2}=\max(\ell_{2},u_{2})
-\max(\ell_{1},u_{1}), X_{3}=\max(\ell_{3},u_{3})-\max(\ell_{2},u_{2}),\ldots$$
are independent and identically distributed with mean $\frac{\sqrt{\pi}}{2}$.
To see this, we have just proved this for the first difference. Now we, so to
speak, move the origin to $(\ell_{1},u_{1})$ and use the homogeneity of the
Poisson process to get that the variables are identically distributed.
Independence follows from the independence properties of the Poisson process.

Therefore, by the Law of Large Numbers, \cite{wiki} for any
$\epsilon>0$ sufficiently small,
\begin{eqnarray*}
\lim_{n\rightarrow\infty}P\left((1-\epsilon)\frac{\sqrt{\pi}}{2}<\frac{X_{1}+\ldots X_{n}}{n}
<(1+\epsilon)\frac{\sqrt{\pi}}{2}\right)=1\\
\Rightarrow \lim_{n\rightarrow\infty}P\left((1-\epsilon)\frac{\sqrt{\pi}}{2}
< \frac{\max(x_{m},y_{m})}{m}<(1+\epsilon)\frac{\sqrt{\pi}}{2}\right)=1
\end{eqnarray*}
just by simplifying the telescoping sum in the definition of the $X_{i}$.

Let $r(n)$ denotes the minimum $r$, such the area $[0,r]^{2}$ contains
exactly $n$ points of the Poisson Process. Then these points determine $n$
random intervals. (Conditional on the number of points of a Poisson process
in a certain area being given, the points themselves are uniformly distributed
over that area).

Recall that we are studying $Y_{n}$, the size of the largest independent
set in a random interval graph with $n$ intervals, that is the longest
chain of intervals in the corresponding interval order. Thus the above
remarks show that we can identify $Y_{n}$ with the largest
$m$ such $(\ell_{m}, u_{m})$ lies in the area $[0, r(n)]^{2}$.

Now, because the Poisson process has density 1, we have
$$\lim_{n\rightarrow\infty}P\left((1-\epsilon)\sqrt{n}<\ell(n)
<(1+\epsilon)\sqrt{n}\right)=1$$
using the Law of Large Numbers \cite{wiki}. Thus, if now we let
\begin{eqnarray*}
m_{1}=[(1-\epsilon)(\frac{2}{\sqrt{\pi}})\sqrt{n}]\mbox{~and~}\\
m_{2}=[(1+\epsilon)(\frac{2}{\sqrt{\pi}})\sqrt{n}]
\end{eqnarray*}
we see by the previous results that, for $n$ sufficiently large, the point
$(\ell_{m1}, u_{m1})$ will lie inside the square $[0, r(n)]^{2}$ and the
point $(\ell_{m2},u_{m2})$ will lie outside of this area, with probability
tending to 1. Thus we indeed get
\begin{eqnarray*}
\lim_{n\rightarrow\infty}P\left(
(1-\epsilon)\frac{2}{\sqrt{\pi}}<\frac{Y_{n}}{\sqrt{n}}
<(1+\epsilon)\frac{2}{\sqrt{\pi}}\right)=1
\end{eqnarray*}
which completes the proof. $\bullet$
\\
\\
We again compare this with the result for $G(n,2/3)$. Here it turns out that
the independence number is, with probability tending to 1 as
$n\rightarrow\infty$, about $2\log_{3}(n)$. Again we refer to
\cite{bollobas} for a proof of this fact. Note again that this number
is much smaller in the $G(n,2/3)$ than in the random interval graph.

\subsection{Comparison of different models}
Here we discuss the differences, which arise in two different models
of random graphs, namely the Erd\H{o}s-R\'{e}nyi model and the
random interval graphs. In the first model, the possibility of an
edge arising is a constant equal to $2/3$, independent from the
number of edges. Moreover from \cite{Schein} we present other
features of random interval graphs, as the value of minimum, maximum
degree for the two models, chromatic and independence numbers. As we
previously saw, for the Erd\H{o}s-R\'{e}nyi model the minimum and
maximum degrees are $\frac{2}{3}n-o(n)$ and $\frac{2}{3}n+o(n)$
respectively and all the degrees are roughly close to $2n/3$. On the
other hand, from \cite{Schein} we have that the minimum degree is
equal to $O(\sqrt{n})$ and the maximum degree is about $n-1$ for the
ordinary model. Moreover chromatic and independence numbers are
$O(\frac{n}{\log n})$ and logarithmic, $O(\log n)$ for
Erd\H{o}s-R\'{e}nyi model, see \cite{Schein}. Also, the clique
number is equal to independence number. Furthermore, in this model
random graphs are not perfect, as the chromatic number is not equal
with the clique number. In random interval graphs chromatic number
and clique number is about to $\frac{n}{2}$, as we saw in above
sections. Finally, from \cite{Schein} independence number is equal
to $O(\sqrt{n})$. This comparison between these models clearly shows
the differences arising in their characteristic values, which
determine their properties.

\pagebreak

\noindent
\section{Variants}

\subsection{Introduction}
In this rather miscellaneous section, we discuss an alternative way
to define random interval graphs, a recent generalization by
Scheinerman of these graphs and some applications of them.

\subsection{Scheinerman's generalization}
Scheinerman has recently introduced a common generalization of both
random interval graphs and the Erd\H{o}s-R\'{e}nyi graphs, namely
random dot product graphs. Also from \cite{scheine} he gives various
definitions of interval graphs.

To understand Scheinerman's idea, we first introduce \emph{Intersection
Graphs}. Suppose we have a finite set of $n$ vertices, $V_{n}$. At each
vertex $v\in V_{n}$, we have a subset $S_{v}\subseteq \mathbf R$ (here,
as usual, $\mathbf R$ is the set of the real numbers). We now say that
two vertices are adjacent if and only if the corresponding sets have a
non-null intersection. In mathematical notation
$$v\sim w\Longleftrightarrow S_{v}\cap S_{w}\neq \emptyset .$$
So an interval graph is a special kind of intersection graph, with
the set $S_{v}$ for each vertex $v$ being an interval of the real line.

Moreover, in \cite{scheine} random intersection graphs are
introduced, by assigning randomly sets $S_{v}$ to the vertices, and then
we again say that two vertices are
adjacent if their corresponding sets intersect. The usual way to assign
these sets is to say that each $S_{v}$ is a subset of $\{1,2,\ldots k\}$
with, for each $v$, $P(i\in S_{v})=p$ say, independently for $1\leq i\leq k$
and each vertex choosing its subset independently. However there are other
possibilities.

A further model studied in \cite{scheine} are \emph{Threshold Graphs}. Here,
For every vertex $v$ , we assign a number $x_{v}$. Then, two vertices
intersect if and only if the sum of the corresponding numbers is $\geq 1$.
Again, in mathematical notation, $v\sim w\Leftrightarrow x_{v}+x_{w}\geq 1$.
Again we can have random threshold graphs by generating the $x_{v}$ in some
random way.

The main business of \cite{scheine} is to give a new definition of a model
of random graphs which combines all these definitions, by using
dot products. Here, vertex $v$ is assigned a $d$-dimensional vector of
real numbers $X_{v}$. Then, two vertices are adjacent, if the corresponding
inner product of the vectors is $\geq 1$. Mathematically,
$$v\sim w\Longleftrightarrow X_{v}\bullet X_{w}\geq 1$$
where $\bullet$ denotes inner (dot) product:
$$(x_{1},x_{2},\ldots x_{d})\bullet (y_{1},y_{2},\ldots y_{d})=\sum_{i=1}^{d}
x_{i}y_{i}.$$

The idea behind the definition of random dot products is that
various ways of defining random interval graphs can be replaced by
the random dot product. Indeed, we have that each vertex $i$ is randomly
assigned a $d$-dimensional vector $X_{i}$. Here $d\in\mathbf N$ is fixed. The
vectors themselves can be generated from some $d$-dimensional distribution:
this could be each component chosen independently, but there are other
possibilities as well. We now say that $i\sim j$ with probability
$f(X_{i}\bullet X_{j})$ for some fixed, and carefully chosen, function $f$.

This is a general definition, which generalizes several of the
definitions above:
\\
\\
Erd\H{o}s-R\'{e}nyi graphs: generalized because if we take, for
every vertex $v$, $X_{v}={\bf x}=(x,x,\ldots x)$ where ${\bf
x}\bullet {\bf x}=p$, and $f(x)=x$, a moment's thought will show
that we recover the $G(n,p)$ model.
\\
\\
Random intersection graphs with each $S_{v}\subseteq \{1,2,\ldots k\}$: because,
if we take $X_{i}$ to be the vector whose $j$th component is $1$ if
$j\in S_{v}$, and whose $j$th component is 0 otherwise, then the property
that two vertices $v$ and $w$ are adjacent if and only if
$S_{v}\cap S_{w}\neq \emptyset$ can be written as the property that
$v\sim w$ with probability $f(X_{v}\bullet X_{w})$, where $f(t)$ is $0$
if $t=0$ and is $1$ otherwise.
\\
\\
(Note: Observant readers will have observed that this is only a
generalization, in the strict sense, of the random intersection
graphs in the case when each $S_{v}\subseteq \{1,2\ldots k\}$,
whereas of course to get intersection graphs to generalize interval
graphs we have to have the $S_{v}$ being infinite sets, namely
certain intervals of the real line. However it is certainly a
generalization in spirit of the idea).

An attractive feature of this very general definition is that we can combine
random and non-random ideas in giving the definitions of the vectors,
according to the situation we are working in.

Scheinerman \cite{scheine} starts by giving some results for the case
when $d=1$ and the \lq vectors\rq\, (really, in this case, scalars, so
we will denote them by the lower case letter) $x_{i}$
are uniformly distributed on $[0,1]$. He takes $f(t)=t^{r}$ for some fixed
$r$: these assumptions will remain in force throughout this paragraph. Now we have
\begin{eqnarray*}
P[i\sim j]=f(x_{i}\bullet x_{j})=\int_{0}^{1} \int_{0}^{1}(x_{i}x_{j})^{r}dx_{i}dx_{j}=\frac{1}{(1+r)^{2}}\\
\end{eqnarray*}
since this is the average, over all possible values of $x_{i}$ and $x_{j}$,
of $f(x_{i}x_{j})$. This is, just by simple integrations,
\begin{eqnarray*}
\int_{x_{i}=0}^{1}x_{i}^{r}dx_{i}\int_{x_{j}=0}^{1}x_{j}^{r}dx_{j}\\
=[\frac{x_{i}^{r+1}}{r+1}]_{x_{i}=0}^{1} [\frac{x_{i}^{r+1}}{r+1}]_{x_{j}=0}^{1}\\
=\left(\frac{1}{r+1}-0\right)\left(\frac{1}{r+1}-0\right)\\
=\frac{1}{(1+r)^{2}}
\end{eqnarray*}
Thus, the expected number of edges is $\frac{n(n-1)}{2}(1+r)^{-2}$,
since it is the expectation of a sum of $n(n-1)/2$ indicator
variables of whether each edge is present, each indicator having
expectation $1/(1+r)^{2}$. He also presents a short calculation, the
details of which we omit, showing that if $a\sim b$ and $b\sim c$,
then conditional on this information $P(a\sim c)$ is larger than it
would be unconditionally: that is, there is a \lq clustering\rq
effect. He believes, but cannot at present prove, that the degrees
in the graph follow a \lq power law\rq: that is, letting $N(d)$
denote the number of vertices of degree $d$, a plot of $\log(N_{d})$
against $\log(d)$ should be a roughly straight line with negative
gradient. Scheinerman observes that various large networks arising
in real life have been observed to have this property (at least
roughly). He further obtains that the expected number of isolated
vertices is $C_{r}n^{(r-1)/r}$ for a suitable constant $C_{r}>0$. In
particular, these graphs are not connected: however, they do have a
very large component and a few isolated vertices. Further they have
diameter at most 6. (The diameter of a graph is the worst case of
the distance between two points in it).

Moreover, Scheinerman \cite{scheine} introduces the \lq inverse
problem\rq . Given a graph on a specific set of vertices, which
vectors are more suitable to model his graph? An obvious approach is
to say that the best choice of $X$s are those which maximize the
likelihood function. This is doable in dimension 1, though in higher
dimensions it becomes very unpleasant fairly rapidly. He thus
suggests an alternative approach based on matrix theory, the Gram
Matrix Approach. In detail: given $G_{1}, G_{2}\ldots G_{m}$ let
$A=\frac{1}{m}\sum_{j=1}^{m}A(G_{j})$
\begin{eqnarray*}
a_{i,j}\approx P[i\sim j]=x_{i}x_{j} (i\neq j)\\
X=[x_{1},x_{2}\ldots,x_{n}]\\
A=X^{t}X\\
\end{eqnarray*}

\subsection{Prisner's definition}
In \cite{prisner}, E. Prisner proposes the following question:
\\
\\
\lq What other reasonable models, apart from Scheinerman's, are there for
random interval graphs? For example, suppose we choose $n$ unit intervals
(i.e. intervals of length 1) which are chosen from the interval $[0,m]$ and
are chosen uniformly and at random?\rq
\\
\\
We are not aware of any substantial work on this model.

\subsection{Applications of Random Interval Graphs}
The obvious applicability of random interval graphs is to scheduling and
assignment problems. For example, suppose each of $n$ people in an
office has an interval each day when he is free for a meeting. The
exact size of this interval, and its position in the period of the working
day (late, or early, or whatever) will vary from day to day, so can
be modeled as random. Then, if we have a random interval graph whose
vertices are the $n$ individuals and where two vertices (individuals) are
adjacent if and only if their random intervals intersect, then we are
saying that these two individuals will have a chance to meet on that day:
and, for example, the largest clique in the random interval graph will be
the largest number of people who can all meet. Similarly, the largest
independent set will be the largest set of people, no two of whom can
meet. Of course, the assumption that the intervals are uniformly distributed
over the working day is probably not very realistic: for example, most
people will be unavailable for some time over lunch, and in practice
they will have several time intervals at which they are available, rather than
just one. (\lq I am available between 9.00 and 10.00, and between 2.00 and 3.00\rq).
However it is a reasonable first model.

Similarly, if we have a series of jobs to carry out in a factory. Suppose
for example we are making a car or similar. Various tasks - say, $n$ of
them - have to be carried out during the production process (for example:
painting the outside, installing the radio, checking the braking system, etc.),
and usually we cannot be doing more than one of these things at once.
There will be time intervals during a working day  in which the people qualified to carry out the
various tasks (brake testing, painting etc.) are available: again, these
will be hard to predict in advance, so can be modeled as random. Again, we
will then want to have a large independent set in the graph, as that means
we have the corresponding time intervals are disjoint, so there are no
clashes.   That is, by doing one of the jobs in its time interval, we
are not reducing our chance of getting one of the other jobs done that day.
(If the painter and the radio installer are available in disjoint time
intervals, then we know that we can just get on with doing the painting
and this will not reduce our chances of getting the radio installed that day
as well).

Moreover, from my essay, \cite{iliop} interval graphs are widely
used in resource allocation problems. That is, we want to allocate a
fixed amount of assets in production activity, in order to maximize
profit. This is a problem in the field of combinatorial
optimization. Another application of interval graphs is their
usefulness in many problems of this field of discrete Mathematics.
An example is the traveling salesman's problem. A salesman leaves
his home and he is willing to visit $n$ towns. He then has to
consider $n!$ alternative feasible tours. We want to find the
optimal tour, so as to visit each town only once and to minimize the
relevant cost of traveling. Obviously, this is a challenging
problem, as the set of feasible tours is too vast. (For example, if
he has to visit 4 towns, the number of feasible tours is 24. If he
has to visit say 6 towns, then we have to consider 6!=720 different
tours!). Consider, to each town we assign a vertex. Then, two
vertices are adjacent if the salesman leaves the town and goes to
the other. Then, to each vertex let's assign an interval. Clearly,
if two intervals intersect, there is an edge arising. But, it is a
cost associated to the salesman tour. So, we aim to find the tour
having minimum cost, so all the $n$ towns will be visited. Also
interval graphs have many other applications, as we shall see.

The textbook by McKee and McMorris \cite{mm} in Chapter 3 contains
references to various applications of interval graphs in Biology,
Psychology and Computing. In Biology, for example, from \cite{mm} a
prominent application of Interval Graphs is the physical mapping of
DNA. From a DNA sequence, some fragments, which are called {\it clones} are
obtained and the goal is to reconstruct the placement of the clones: that is,
where they are on the DNA string. Thus a clone is an interval of a line of
DNA. 

To turn this into a problem about random interval graphs, we assign to
each clone a vertex. Two different clones are
adjacent if and only if their corresponding intervals intersect. This
clearly gives an interval graph.

Also, interval graphs are used in social sciences. For example, in
Psychology, \cite{mm} they are widely used as tools, measuring
notions, which determine different psychological theories. (Most
theory of measurement is based on physical science: however, in the
social sciences, different theories may be more appropriate).

An example from \cite{mm} is that a person has a set $A$ of
alternatives solutions to choose. For simplicity, suppose that
elements of $A$ are different makes of cars. Our person has
preferences among the the different makes of cars: for example,
he might prefer stylish cars or cars having low cost of service, etc). Then,
a real-valued function $f$ on the set $A$, such that for $a,b\in A$,
he prefers alternative $a$ than $b$ when $f(a)>f(b)+\delta$, where
$\delta$ is a positive constant representing a threshold - a \lq just
noticeable\rq\, difference between the two kinds of cars. Then, we
define a binary relation $R$ on $A$ to be an
interval order on the set $A$ of alternatives, if it satisfies the
two following axioms. (We are thinking here of $aRb$ as meaning that
$a$ is preferable to $b$). 
\\
\\
\begin{Large}
Axiom 1:
\end{Large} For all $a\in A$, not $aRa$
\\
\begin{Large}
Axiom 2:
\end{Large}For all $a,b,c,d\in A$, if $aRb$ and $cRd$, then either $aRd$ or $cRb$\\
\\
The motivation is that a car is not preferable to itself (clearly), which
gives Axiom 1. Similarly, if $a$ is preferable to $b$ and $c$ preferable to
$d$, it seems reasonable that at least one of $a$ is preferable to $d$ and
$c$ is preferable to $c$ should hold. 

Furthermore, interval graphs are used in Computing. They are used in
scheduling problems. From \cite{mm} we have an interesting
application of this class of problems. Suppose that we want to find
an arrangement, in order to construct a timetable for different
courses in a University. We have a fixed number of rooms available
for teaching purposes and we know the number of teachers. We aim to
construct an efficient timetable, so to be no overlap between
teaching hours for every lesson. We assign various courses to
vertices. Then, two vertices are adjacent, when the corresponding
intervals intersect. When this is the case, we have two different
courses at the same hour. Thus, we want to find the minimum number
of rooms needed, in order all the courses to be taught. In "graph
language" we want to find the chromatic number. That is the minimum
number of colors needed, so two connected vertices have different
colors, \cite{mm}. Moreover, in problems related to information
retrieval, we use interval graphs. Suppose that $\Phi$ denotes set
of files, which contain information and $Q$ is the set of queries
for retrieving information. Then, $\Phi$ and $Q$ satisfy the
\emph{consecutive retrieval property} if the files relevant to each
query can be stored consecutively in a linear form, so not to be
overlap, \cite{mm}. We easily deduce from the above the extensive
use of interval graphs in various, different sciences, from areas of
applied Mathematics to social sciences, as Psychology.

\pagebreak

\noindent
\section{Conclusions: areas for further work}
In this project, we examined Random Interval Graphs, emphasizing
Scheinerman's definition. We presented some results on the number of
edges and then considered the degrees of vertices in the graphs,
observing that these are much more spread out that in the
alternative Erd\H{o}s-R\'{e}nyi model of random graphs: for example,
the maximum degree of a random interval graph is very likely to be
close to $n-1$, and indeed is equal to $n-1$ with probability $2/3$.
We then considered cliques, independent sets and chromatic numbers
of random interval graphs, obtaining asymptotic estimates of each of
the quantities involved and comparing their values with the values
in the Erd\H{o}s-R\'{e}nyi model. Finally, we wrote about other ways
of defining Random Interval Graphs, such as Scheinerman's random dot
product graphs, emphasizing the 1-dimensional case of Scheinerman's
theory. Finally we outlined some areas of application.

There are some questions left unanswered by our work. For example,
it would be desirable to investigate measures of connectivity (such
as vertex-connectivity or edge-connectivity) in random interval
graphs. Also: what is the diameter of a random interval graph? The
obvious guess would be that it is 2, since if we have two vertices
$v$ and $w$, one would hope that one of the many vertices of high
degree (close to $n-1$) will be adjacent to both of them. Certainly
the probability that the diameter is 2 is at least 2/3, since if
there is a vertex of degree $n-1$ is till be adjacent to both $v$
and $w$ (if either $v$ or $w$ is a vertex of degree $n-1$, then it
is clearly adjacent to the other). The result that the probability
is at least 2/3 is now just a consequence of the fact that with
probability 2/3 there is a vertex of degree $n-1$. However the guess
that the diameter of a random interval graph is 2 with probability
tending to 1 does not seem to follow immediately from what we have
proven, because it is not quite clear that we can avoid the
situation where there are two vertices $v$ and $w$ of low degree and
all the vertices of high degree fail to be adjacent to at least one
of $v$ and $w$.

Another topic is the existence of Hamilton cycles in a random
interval graph. Scheinerman (\cite{Schein}) shows, by a rather long
and difficult argument, that with probability tending to 1, a random
interval graph is Hamiltonian: that is, it has a cycle which passes
through every vertex of the graph. A $G(n,2/3)$ is also Hamiltonian:
indeed $G(n,2/3)$ has the stronger property that it has (with
probability tending to 1) $\lfloor \delta(G)/2\rfloor$ edge-disjoint
Hamilton cycles: we refer to \cite{bollobas} for a proof of this
fact. (Two cycles are edge-disjoint if and only if there is no edge
which is in both cycles). Note that $\lfloor \delta(G)/2\rfloor$ is
the largest number of edge disjoint Hamilton cycles we could have in
a graph $G$, because each Hamilton cycle will use up two edges in
passing through a vertex $v$ of degree $\delta(G)$. Is it true that
a random interval graph has $\lfloor \delta(G)/2\rfloor$
edge-disjoint Hamilton cycles? Again this does not seem to be
obvious.

We hope that we have given a reasonable material of random interval
graphs, some interesting results on them, and  some ideas of how
they might be useful.

\pagebreak

\end{document}